\documentclass[letterpaper, 10 pt, conference]{ieeeconf}
\IEEEoverridecommandlockouts
\usepackage{cite}
\usepackage{amsmath,amssymb,amsfonts}
\usepackage{algorithmic}
\usepackage{graphicx}
\usepackage{textcomp}
\usepackage{xcolor}
\usepackage{dsfont}
\usepackage{subfig}
\usepackage{makecell}

\usepackage{scalefnt}
\usepackage{tikzscale}
\usepackage{tikz}
\usepackage{pgfplots}
\pgfplotsset{compat=1.17}
\usepackage{changepage}
\usetikzlibrary{shapes.geometric}
\usetikzlibrary{calc}
\usetikzlibrary{positioning}
\usetikzlibrary{intersections}
\usetikzlibrary{decorations.fractals,spy}

\def\BibTeX{{\rm B\kern-.05em{\sc i\kern-.025em b}\kern-.08em
    T\kern-.1667em\lower.7ex\hbox{E}\kern-.125emX}}
\begin{document}

\title{Safety Envelope for Orthogonal Collocation Methods in Embedded Optimal Control
}



\author{{Jean Pierre Allamaa$\,\,^\textrm{1,2}$, Panagiotis Patrinos$\,\,^\textrm{2}$, Herman Van der Auweraer$\,\,^\textrm{1}$, Tong Duy Son$\,\,^\textrm{1}$}
	\thanks{This project has received funding from the European Union’s Horizon 2020 research and innovation programme under the Marie Skłodowska-Curie grant agreement ELO-X No 953348}
	\thanks{$^\textrm{1}$  Siemens Digital Industries Software,  3001 Leuven, Belgium. Email: \tt jean.pierre.allamaa@siemens.com}
	\thanks{$^\textrm{2}$  Dept. Electr. Eng. (ESAT) - STADIUS research group, KU Leuven, 3001 Leuven, Belgium. Email: \tt panos.patrinos@esat.kuleuven.be}
}
\maketitle

\begin{abstract}
Orthogonal collocation methods are direct approaches for solving optimal control problems (OCP). A high solution accuracy is achieved with few optimization variables, making it more favorable for embedded and real-time NMPC applications. However, collocation approaches lack a guarantee about the safety of the resulting trajectory as inequality constraints are only set on a finite number of collocation points. In this paper we propose a method to efficiently create a convex safety envelope containing the trajectory such that the solution fully satisfies the OCP constraints. We make use of Bernstein approximations of a polynomial’s extrema and span the solution over an orthogonal basis using Legendre polynomials. The tightness of the safety envelope estimation, high accuracy in solving the underlying differential equations, fast rate of convergence and little conservatism are properties of the presented approach making it a suitable method for safe real-time NMPC deployment. We show that our method has comparable computational performance to pseudospectral approaches and can accurately approximate the original OCP up to 9 times more quickly than standard multiple-shooting method in autonomous driving applications, without adding complexity to the formulation.
\end{abstract}

\section{Introduction}
Nonlinear model predictive control (NMPC) for autonomous systems has been extensively studied in academia and industry as an advanced control strategy incorporating by design safety and comfort, and capable of achieving high-performance. NMPC has shown promising results in engineering applications such as robotics and autonomous driving (AD), specially in the context of safety critical scenarios where constraints have to be satisfied alongside an optimal decision making~\cite{9867514}. However, the main challenge remains in the transfer of the NMPC to the real-world with critical real-time (RT) implementation. 


At its core, NMPC requires the solution of an optimal control problem (OCP) at every iteration. This entails solving a boundary value problem, satisfying state and input constraints, and minimizing an objective function. Most dynamical systems in MPC, are first formulated using ordinary differential equations (ODE) in continuous time, to capture the physical system dynamics. The performance, accuracy and RT capabilities depend on the choice of the ODE's discretization scheme, which is often neglected. Rather, the current focus is on the solver and optimization routine approaches that can be solved using one of the available methods in the literature~\cite{10.5555/1942009, doi:10.1137/1.9780898718577, rawlings2017model}. For embedded applications, some methods are proposed such as the real-time iteration (RTI) scheme~\cite{DIEHL2002577}, single shooting for low-dimensional systems~\cite{8550253}, and explicit MPC. In most applications, the continuous problem is approximated by its discrete time counterpart using a simple, low accuracy integrator~\cite{xue_active_2020}, or more accurate schemes such as the Runge-Kutta 4\textsuperscript{th} order method (RK4) which are not always suitable for RT control.
We propose a collocation approach for transcribing an OCP with fewer variables and without discretization, benefiting from a fast rate of convergence and highly accurate smooth solution suitable for RT applications, such as in AD. Albeit using fewer variables, our method guarantees a safe trajectory without additional computational burden or major conservatism.

Collocation is a method for solving non-linear OCPs~\cite{467672, Li2016PseudospectralOC}, and ODEs more efficiently than traditional finite elements as the discretization schemes are omitted and the computational effort of explicitly solving ODEs at each iteration is saved~\cite{BIEGLER1984243}. It avoids discretization by approximating the solution with splines or polynomials, infinitely differentiable, continuous, suitable for smooth dynamical systems. This is crucial for autonomous systems with safety critical operation, as it contributes to the stability, safety and performance. The most common collocation methods belong to the family of pseudospectral or orthogonal collocation with fast rate of convergence~\cite{BIEGLER1984243}. However, the main problem remains to this date that transcription from OCP to nonlinear programming problem (NLP) in the collocation approach only sets the inequality constraints at the collocation nodes~\cite{garg:hal-01615132, doi:10.1080/00207177508922030}. Thus, the safety of the continuous trajectory is not guaranteed. To overcome this issue, some researchers introduce slack variables on the inequality constraints~\cite{rabiei_collocation_2020, Raghunathan2004DynamicOS}. Others use B-splines with the convex hull property, which introduces conservatism based on the distance between the control polygon and the spline~\cite{8022960}. 

In order to address those challenges and enhance safety without a compromise on computation time and solution accuracy, the contributions of this paper are: 1) a collocation approach aimed towards embedded application resulting in a speed up in the numerical optimization scheme allowing to deploy NMPC in real-time, orders faster than standard shooting methods through 2) an orthogonal collocation (SOC) with safety envelope SOCSE for a tight estimation of the spline's extrema at optimization runtime allowing to constrain the full trajectory within the OCP boundaries with the least conservatism and guaranteeing a feasible everywhere optimal trajectory 3) a continuous spline solution allowing to span a multitude of trajectories, benefiting from fast convergence speed and high accuracy in solving the ODE and OCP.   

This paper is organized as follows: Section~\ref{section:2} discusses a background on the different methods to solve an OCP. Section~\ref{section:3} introduces the simultaneous collocation-based approach for solving ODEs. Furthermore, section~\ref{section:4} presents the proposed method for generating a convex safety envelope for the trajectories in spectral orthogonal collocation and validates the approach in simulation on an academic example. Finally, the SOCSE approach is employed in a complex autonomous driving application in simulation and the method is deployed in an embedded vehicle-hardware-in-the-loop autonomous valet parking application with the results presented in section~\ref{section:5}.
\section{Optimal control problem}
\label{section:2}
Optimal control problems are initially posed in continuous time and seek to optimize a cost function over the state $x$ and control actions $u$, while satisfying a set of (nonlinear) inequality and equality constraints. 
For the remainder of the work, the focus is on the nonlinear continuous Bolza problem~\cite{Huntington2007AdvancementAA}:
\begin{equation}
	\label{eq:generic_ocp}
	\left\{
	\begin{aligned}
		\min_{x(.),u(.)} &J(x,u) = \phi(x(t_f)) + \int_{t_0}^{t_f} l(x(t),u(t))dt \\
		\textrm{subject to } &\dot{x}(t) = f(x(t), u(t)),  \\
		&g(x(t), u(t)) \leq 0,\\
		&g_{f}(x(t_f)) \leq 0, \\
		&x(t_0) = \bar{x}.
	\end{aligned}
	\right.
\end{equation}
This infinite dimensional generic OCP in~\eqref{eq:generic_ocp} is transcribed into a finite dimensional NLP, most often with discrete time variables on a time grid in $[t_0,t_f]$ of size $N_H$ steps. The dynamics $f(x,u)$ are discretized into $f_d$ and the integral of the Lagrange term $l(x,u)$ is approximated by a finite sum over the discrete nodes.
Moreover, $\phi$ is the terminal cost, $g(x,u)$ are the nonlinear constraints and $g_{f}$ is the terminal constraint.
The accuracy of the solution of the discrete OCP depends on the number of elements $N_H$ in the horizon. Usually $N_H$ is derived as $(t_f-t_0)/T_s$ where $T_s$ is the sampling time. Therefore, to truly benefit from the predictive capabilities of MPC, $N_H$ is large. Albeit its accuracy, this comes with a computational burden and limits the RT capability, which we aim to tackle. In this chapter, the three most common approaches for solving a generic NLP are presented.
\subsection{Discrete-time conventional approaches}
Direct methods optimize over the discretized OCP, and can deal with non-smooth systems by dividing the time domain into smaller subintervals.
\paragraph{Single shooting, a sequential approach}
Single shooting (SS) alleviates the computational burden by eliminating the state sequence $(x_0, \dots, x_{N_H})$ and expressing the minimization in terms of the control actions $(u_0, \dots, u_{N_H-1})$ and the state trajectory is rolled out using $f_d$~\cite{8550253}. The NLP can be efficiently solved with solvers exploiting this structure such as \textit{PANOC}~\cite{8263933}. 
\paragraph{Multiple shooting, a hybrid simultaneous approach}
SS does not scale well with the dimensionality of the system as for complex systems, the sparsity structure in the Hessian is lost. Another direct approach, is Bock's multiple-shooting method~\cite{Bock1984AMS}. It solves the discrete time OCP simultaneously for $x$ and $u$, and benefits from structure and sparsity exploitation. 
Similar to SS, the NLP can be solved using a Sequential Quadratic Programming (SQP) method or an interior point method solver, such as IPOPT. For embedded applications, SQP is usually a more suited method, such as in the RT autonomous driving applications in~\cite{9867514, 8061009}. 
\subsection{Collocation: Simultaneous approach}
A third approach for solving the OCP, is the collocation method.
It is a fully simultaneous direct approach which transcribes the continuous time OCP in~\eqref{eq:generic_ocp} into an NLP without numerical integration schemes or embedded ODE solver. Instead, the trajectory is assumed to be of polynomial form and the generic solution is of the form:%
\begin{equation}%
	\label{eq:generic_solution_form}
	x(t) = \sum_{k=0}^{K}\alpha_k\phi_k(t),
\end{equation}%
where $\alpha_k$ are the expansion coefficients and $\phi_k(t)$ is a spanning basis.
In pseudospectral collocation (PS), the state and control trajectories are represented by one finite element polynomial of degree $M$ and the solution is uniquely determined given the initial condition $x_0$. That is, the NLP solves for the state value $\alpha_k = x(\tau_k)$ at a collocation node $\tau_k$ and $\phi_k(\tau)$ is most commonly a Lagrange interpolating polynomial. The advantage of this method lies in its use of a spectral grid, which converges quickly and with high accuracy to the true solution of the ODE, assuming smooth trajectories~\cite{rawlings2017model}. It is more favorable than a uniform grid that has undesirable properties in polynomial approximation~\cite{Huntington2007AdvancementAA}. 
\section{Spectral and orthogonal collocation methods for solving ODEs}
\label{section:3}
In this chapter, we present an approach to solve the OCP via spectral orthogonal collocation, by first defining the components in the OCP to NLP transcription: spectral grid forming, Gaussian quadrature and orthogonal collocation.

In a collocation approach, the problem is a dynamic optimization with polynomial representation of the state and control in each finite element~\cite{10.5555/1942009}. In particular, Legendre polynomial based collocation will be used in the remainder of this work. Under this design, truncation errors for solving the ODEs and fitting a polynomial over the collocation points' data are small~\cite{10.5555/1942009}. However, a drawback in collocation methods, is that inequalities are usually only set at the collocation points, with no guarantee that the resulting polynomial is continuously within the OCP boundaries. We propose to tackle this by setting the solution as a spline and efficiently shaping the trajectories through the spline's convex safety envelope.

\subsection{Spectral collocation and Gauss quadrature}
\label{subsection:quadrature_spectral}
In spectral-based collocation approaches, an accurate solution is obtained using few collocation nodes with a spectral accuracy of convergence, faster than any power of $1/N$ for smooth problems~\cite{Huntington2007AdvancementAA}, with $N$ the number of collocation points. There are three families of Legendre-Gauss quadratures to be used in a spectral collocation, among which the Legendre-Gauss-Lobatto (LGL). For this work, the focus is on LGL as it contains both ends $t_0$ and $t_f$, but the presented approach generalizes beyond the choice of quadrature.
Readers are referred to other variants in~\cite{garg:hal-01615132}. The collocation points $\tau_i, i=2,\dots,N-1$ are uniquely chosen to be the roots of $\dot{\mathcal{L}}_{N-1}$, the derivative of the Legendre polynomial $\mathcal{L}_{N-1}$ of degree $N-1$, c.f~\eqref{eq:legendre_polynomial}. The two remaining collocation points are the boundaries $[-1,1]$ such that $\tau_i =: \begin{bmatrix} -1 & \textrm{Roots of } \dot{\mathcal{L}}_{N-1} & 1\end{bmatrix}$~\cite{doi:10.1137/1.9780898718577}. This unique choice of collocation points defines the spectral grid of a Lobatto PS method (LPM)~\cite{garg:hal-01615132}. LPM can achieve exact estimation up to machine precision, as long as the true polynomial is up to degree $2N-3$, with a relatively small number of nodes~\cite{Huntington2007AdvancementAA} and a truncation error $O(h^{2N-2})$~\cite{10.5555/1942009}. Since the Legendre polynomial is defined over $\tau\in[-1,1]$, the transformation from the original time domain $t\in[t_0, t_f]$ is: $t=\frac{t_f-t_0}{2}\tau + \frac{t_f+t_0}{2}$. %
It follows for scaling the dynamics and cost integral in $\tau$: $\int_{0}^{t_f}l(x(t),u(t))dt = \frac{t_f}{2}\int_{-1}^{1}l(x(\tau),u(\tau))d\tau$ with the time scale $t_f/2$. %
Finally, the integral of the Lagrange term $l(x(\tau),u(\tau))=:l(\tau)$ in the OCP cost function is given by the exactness of the Gauss quadrature rule as:%
\begin{equation}%
	\label{eq:gauss_quadrature}
	\int_{-1}^{1}l(\tau)d\tau = w_1 l(-1) + w_N l(1) + \sum_{i=2}^{N-1}w_il(\tau_i),
\end{equation}%
where $w_i$ are the LGL quadrature coefficients. They are computed as in~\cite{Huntington2007AdvancementAA}:%
\begin{align}%
	\label{eq:gauss_quadrature_coefficients}
	w_1 &= w_N = 2/(N(N-1)), \nonumber \\ 
	w_i &= \frac{2}{N(N-1)\mathcal{L}_{N-1}(\tau_i)}, i=2,\dots,N-1,
\end{align}%
The accuracy of the Gaussian quadrature in~\eqref{eq:gauss_quadrature} can provide an exact solution to the integral in~\eqref{eq:generic_ocp}~\cite[Theorem 10.1]{10.5555/1942009} and is valid for the presented choice of spectral grid $\tau_i$.
\subsection{Orthogonal collocation}
The common PS formulation using a single element Lagrange polynomial gives little flexibility in terms of shaping the path beyond the collocation points. For this, we exploit the spectral property and the Gauss quadrature in~\eqref{eq:gauss_quadrature} and~\eqref{eq:gauss_quadrature_coefficients} in a orthogonal collocation framework with constrained splines. 
Instead of limiting the solution to one polynomial, it is possible to set the trajectories as piece-wise combination of polynomials. Hence, the solution, and its derivatives, are smooth and continuous by construction. For orthogonality, we propose a solution of the form in~\eqref{eq:generic_solution_form} spanned by a basis of Legendre polynomials $\phi_k(\tau) = \mathcal{L}_k(\tau)$. Legendre polynomials $\mathcal{L}_k$ are orthogonal, satisfy $\mathcal{L}_k(1) = 1$ for any degree $k$ and are given by the Rodrigues formula~\cite{Huntington2007AdvancementAA}:%
\begin{equation}%
	\label{eq:legendre_polynomial}
	\mathcal{L}_k(\tau) = \frac{1}{2^k k!} \frac{d^k}{d\tau^k}[(\tau^2-1)^k].
\end{equation}%
We define the orthogonal collocation scheme based on the Legendre-spline to approximate the solutions $\tilde{x}(\tau)$ and $\tilde{u}(\tau)$:%
\begin{equation}%
	\label{eq:legendre_spline}
	\tilde{x}(\tau) = \sum_{k=0}^{M} \alpha_k\mathcal{L}_k(\tau) = \alpha^\top\mathbf{L}_Mv(\tau), (\tau\in[-1,1]),
\end{equation}%
where $\mathbf{L}_M$ is an upper trianglular matrix $\in \mathbb{R}^{(M+1)\times(M+1)}$ formed by the coefficients of $\mathcal{L}_k$ with respect to time $\tau$. $N_x, N_u$ are the number of state and control variables respectively, $M$ is the degree of the fitting spline, $N$ is the number of collocation points. Moreover, $\alpha = \begin{bmatrix}\alpha_0 & \cdots & \alpha_M\end{bmatrix}^\top \in \mathbb{R}^{(M+1)\times N_x}$ and $v(\tau) = \begin{bmatrix} 1 & \tau & \tau^2 & \cdots & \tau^M \end{bmatrix}^\top$ is the time vector. For the remainder of the work,~\eqref{eq:legendre_spline} is referred to as the Legendre-spline of degree $M$. 
We set the $N$ collocation points $\tau_i$ as the LGL nodes presented in \ref{subsection:quadrature_spectral} and differentiate the spline in~\eqref{eq:legendre_spline} at the collocation points $\tau_i$ to satisfy the dynamics:
\begin{equation}
	\dot{\tilde{x}}(\tau_i) = \alpha_x^\top\mathbf{L}_M\dot{v}(\tau_i) = \frac{t_f}{2}f(x(\tau_i), u(\tau_i)).
\end{equation}
Under this formulation, there are more degrees of freedom as the number of collocation points $N$ is not pre-defined for a certain spline of degree $M$. The only requirement to avoid over constraining the NLP is $(N_u + N_x)\times (M+1) \geq N_x\times (N+1)$. %
The right hand side of this inequality is formed by the initial condition and the dynamics equality constraints at the $N$ collocation points. A spline of degree $M$ has $M+1$ coefficients, leading to the left hand side of the equation. 
\section{Safety envelope for collocation}
\label{section:4}
The transformation to a Legendre-spline collocation allows to efficiently approximate the path's extrema for any spline degree. This permits to exploit the spline's bounds inside the NLP without analytical derivation of the maximum and minimum at runtime. The safety envelope presented in this chapter contributes to a numerical speed-up while guaranteeing a safe continuous solution, without additional formulation complexity.
We demonstrate the proposed approach with a simple analytical example showing the effectiveness compared to the shooting methods, and in a later stage on a complex autonomous valet parking example. 
\subsection{Safety envelope through polynomial bounds}
Let $P(t)=a_0+a_1t+a_2t^2+\dots+a_Mt^M$ be a polynomial of degree $M\geq0$ with real coefficients $a_j, j=0,\dots\,M$. Then it follows from the Bernstein form of a polynomial in the work of~\cite{CargoBernstein} that the true minimum $\underline{P}$ and maximum $\overline{P}$ of the polynomial for $t\in[0,1]$ are bounded as:%
\begin{align}%
		\label{eq:bernstein_extremum}
		&\min \{b_j\} \leq \underline{P} \leq P(t) \leq \overline{P} \leq \max \{b_j\}, 
		\\ &\qquad 0\leq t\leq1, j=0,\dots,M.\nonumber
\end{align}%
The Bernstein extrema estimations $b_j$ are defined as:%
\begin{subequations}
	\label{eq:bernstein_coefficients}
	\begin{align}
		b_j &= \sum_{k=0}^{M}a_j \binom{j}{k}\bigg/\binom{M}{k}, \, j=0,1,\dots,M,\\
		b &= \mathcal{B}a, \, \mathcal{B} \in \mathbb{R}^{(M+1) \times (M+1)}.
	\end{align}
\end{subequations}
$\min \{b_j\}$ and $\max \{b_j\}$ define the convex envelope containing the polynomial. Note that if $\min \{b_j\}=b_o$ and $\max \{b_j\}=b_M$, then the approximation of the extrema is exact, with $\underline{P}=b_0$ and $\overline{P}=b_M$~\cite{CargoBernstein}.
Equations~\eqref{eq:bernstein_extremum} and~\eqref{eq:bernstein_coefficients} are defined in the time interval $\tau_* \in [0,1]$ whereas the Legendre-spline is defined in $\tau \in [-1,1]$, with $\tau=(2\tau_*-1)$.
We convert the vector $v(\tau)$ in~\eqref{eq:legendre_spline} to $[0,1]$ using the Binomial theorem:
\begin{subequations}
	\begin{align}
		\label{eq:binomial_time_transform_a}
	    (2\tau_*-1)^M = \sum_{r=0}^{M}\binom{M}{r}(-1)^{M-r}(2)^r\tau_*^r	= \mathcal{E} v(\tau_*),\\
   		\label{eq:binomial_time_transform_b}
	   	v(\tau) = [1,\dots,\tau^M]^\top = \mathcal{E} v(\tau_*), \, \mathcal{E} \in \mathbb{R}^{(M+1) \times (M+1)},
	\end{align}
\end{subequations}
where $\mathcal{E}$ is upper trianglular. Substituting~\eqref{eq:binomial_time_transform_b} in~\eqref{eq:legendre_spline}:
\begin{equation}
	\tilde{x}(\tau_*)=\alpha^\top\mathbf{L}_M\mathcal{E}v(\tau_*).
\end{equation}
The transformed coefficients $\alpha_*^\top = \alpha^\top\mathbf{L}_M\mathcal{E}$ are used in~\eqref{eq:bernstein_coefficients}. Finally, the resulting Legendre-spline is contained in its convex envelope, having as its boundaries the maximum and minimum elements of $\mathcal{P}_M$:
\begin{align}
	\label{eq:safety_envelope}
	\mathcal{P}_M = \mathcal{B}\alpha_* = \mathcal{B}\mathcal{E}^\top\mathbf{L}_M^\top\alpha = \mathcal{C}_M\alpha.
\end{align}
The state and control inequality constraints in the OCP are transcribed into the NLP by adding the $M+1$ linear constraints on the coefficients $\alpha$ in~\eqref{eq:safety_envelope}, guaranteeing the feasibility of the trajectory everywhere in the interval.
The benefit in this method is that the matrix $\mathcal{C}_M = \mathcal{B}\mathcal{E}^\top\mathbf{L}_M^\top \in\mathbb{R}^{(M+1)\times (M+1)}$ is computed at construction time once, given the spline degree $M$, and replaces the state and input constraints at the collocation points in the optimization framework, hence the claim of no added formulation complexity. The resulting NLP is:%
\begin{equation}%
	\label{eq:SOSCSE_nlp}
	\left\{
	\begin{aligned}
		\min_{\alpha_x, \alpha_u} &\frac{t_f}{2}[\sum_{i=2}^{N-1}w_il(\tau_i)+ \phi(\alpha_x^\top\mathbf{L}_M v(1))\\
		&\qquad + w_1 l(x_0, \alpha_u^\top\mathbf{L}_M v(-1)) + w_N l(1)]\\
		\textrm{s.t. } &\alpha_x^\top\mathbf{L}_M v(-1) = \bar{x}_0,\\
		&\alpha_x^\top\mathbf{L}_M\dot{v}(\tau_i) = \frac{t_f}{2}f(\alpha_x^\top\mathbf{L}_M v(\tau_i), \alpha_u^\top\mathbf{L}_M v(\tau_i))\\ 
		&\underline{x} \leq (\mathcal{C}_M \alpha_x)^\top\leq \overline{x},\,\quad \alpha_x\in\mathbb{R}^{(M+1)\times N_x},\\
		&\underline{u} \leq (\mathcal{C}_M \alpha_u)^\top\leq \overline{u},\,\quad \alpha_u\in\mathbb{R}^{(M+1)\times N_u},\\
		&g_{f}(\alpha_x^\top\mathbf{L}_M v(1)) \leq 0,
	\end{aligned}
	\right.
\end{equation}%
where $l(\tau_i) = l(\alpha_x^\top\mathbf{L}_M v(\tau_i), \alpha_u^\top\mathbf{L}_M v(\tau_i)), \forall i \in [1,N]$.
Therefore, by formulating the problem as in~\eqref{eq:SOSCSE_nlp}, the trajectories within the optimization routine, are defined by their extrema, as a function of their coefficients. In addition, the tightness of the extremum approximation using the Bernstein form in~\cite{CargoBernstein}, permits to fit exactly, the complete Legendre-spline trajectory within the inequality bounds of the original OCP formulation, with minimal conservatism. 

\begin{figure}%
	\begin{center}
		{\input{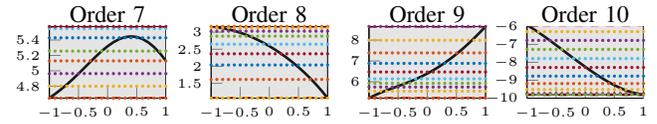}}    
		\caption{Convex envelope of random Legendre-splines with varying degrees. Solid black: Legendre-spline, Grey area: convex envelope estimation, Dotted: $b_j$ elements from~\eqref{eq:bernstein_coefficients}
		\label{fig:envelope_random_splines}}
	\end{center}
\end{figure}

Furthermore, to visualize the tightness of approximation of the safety envelope, we generate random Legendre-splines as in~\eqref{eq:legendre_spline} and form their convex safety envelopes. The results are shown in Figure~\ref{fig:envelope_random_splines}. The elements of~\eqref{eq:safety_envelope} are dotted and the shaded grey area shows the convex set containing the spline. For all the examples except for degree 7, the convex set estimates with no error the true spline's extrema. Note that~\eqref{eq:SOSCSE_nlp} uses one spline per horizon, but it is possible to divide the horizon into multiple segments for non-smooth problems.

\subsection{Analytical academic example}
\label{subsect:analyt_ex}
The following metrics are used for the solution accuracy:
\begin{itemize}
	\item the true analytical solution $\{x^*(.), u^*(.)\}$, and cost $J^*$
	\item Optimal cost deviation =: $100\times(J_{sol} - J^*)/J^*$
	\item Optimal control RMSE =: $\lVert u_{sol} - u^* \rVert^2_2/N^*$, where $u_{sol}$ is the solution control vector sampled at $N^*$ steps
\end{itemize}%
Consider the classical optimization example presented in~\cite{1100470} with the method of weighted residuals, and augmented with constraints in~\cite{doi:10.1080/00207177508922030}. The authors demonstrate the accuracy of the collocation method in solving the OCP, with a small relative error compared to the analytical solution. However, both examples use a non-orthogonal basis, and the choice of collocation points is random. In this work, to demonstrate the effectiveness of the SOCSE method, we propose a modified constrained problem with the input constraint $-0.3 \leq u(t) \leq -0.1\,\forall t\in[0,1]$, such that the optimal solution is affected:%
\begin{align}%
	\label{eq:academic_example}
	\min_{x(.),u(.)} &J(x,u) = \frac{1}{2}\int_{0}^{1} [x^2(t)+u^2(t)]dt \nonumber \\
	\textrm{subject to } &\dot{x}(t) = -x(t) + u(t), \;\forall t\in [0,1]\\
	&0.2 \leq x(t) \leq 1.0\nonumber\\
	-&0.3 \leq u(t) \leq -0.1,\nonumber\\
	&x(0) = 1.0. \nonumber 
\end{align}%

\begin{figure}%
	\begin{center}
		{\input{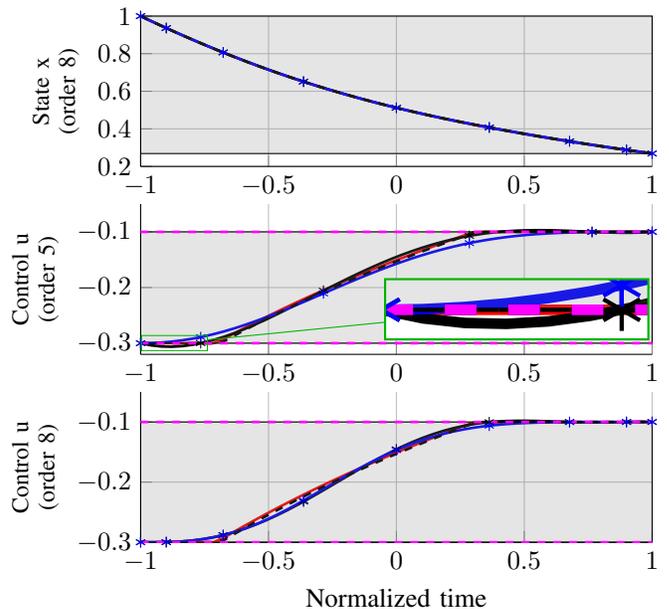}}    
		\caption{Modified academic example: Optimal state and control using different methods (Red: multiple-shooting with 100 nodes, Solid black: Pseudospectral, Blue: SOCSE, Grey: convex envelope, Dashed black: analytical solution, Pink: OCP bounds). PS violates constraints between collocation points.}
		\label{fig:dummy_ex_states_control}
	\end{center}
\end{figure}

\begin{table}[htbp]%
	\caption{Comparison of different methods for solving an academic OCP in~\eqref{eq:academic_example}}
	\begin{center}
		{\renewcommand{\arraystretch}{1.1}
		\scalebox{0.8}{
		\begin{tabular}{p{0.2\linewidth} p{0.2\linewidth} p{0.2\linewidth} p{0.2\linewidth} p{0.2\linewidth}}
			\hline
			Method &Solve time (s)& Optimal cost deviation (\%)& Optimal control RMSE & Continuous OCP bounds violation \\
			\hline
			MS-50 steps& 0.043&2.909 & 2.0e-3 & No\\
			\hline
			PS $\mathcal{O}$(5)& 0.016 & 2.0e-3 &4.0e-3  & Yes \\
			\hline
			PS $\mathcal{O}$(8)& 0.035&0.019 & 4.7e-3  & Yes\\
			\hline
			\textbf{SOCSE $\mathcal{O}$(5)}& \textbf{0.016}&\textbf{0.049} & \textbf{7.4e-3} & No \\
			\hline
			\textbf{SOCSE $\mathcal{O}$(8)}& \textbf{0.024}&\textbf{0.024} &\textbf{4.2e-3}  & No\\
			\hline
		\end{tabular}}}%
		\label{tab:toy_example}
	\end{center}%
\end{table}%
The NLP is formulated as in~\eqref{eq:SOSCSE_nlp} and is solved to convergence using IPOPT with a feasibility tolerance of $1\text{e-}9$. The comparison results are shown in Table~\ref{tab:toy_example} and the analytical solution is derived from Pontryagin's Minimum Principle theorem. State $x(\tau)$ and input $u(\tau)$ are shown in Figure~\ref{fig:dummy_ex_states_control}. For comparison, we solve the problem with our proposed approach (SOCSE), a multiple-shooting (MS) approach with 100 steps, and a pseudospectral collocation (PS) with Lagrange polynomial. The problem is solved with a cold-started solution equal to zero. Both collocation approaches are solved with a degree 5 ($\mathcal{O}$(5), $M=5, N=6$) and 8 ($\mathcal{O}$(8), $M=8, N=9$) spline/polynomials. The optimal state trajectory is indistinctive between all three methods. With degree 5, the PS approach constraining only the collocation nodes fails to provide an over-all safety of the trajectory as seen between time -1 and -0.5. A control action sampled from this interval is therefore unsafe.  However, in our approach, the spline's safety envelope tightly estimates its extrema (grey box) and sets them at the limits of the OCP inequalities (pink lines), over the complete time domain. Using SOCSE $\mathcal{O}$(8), and in comparison with the analytical solution, an accurate OCP solution is achieved, with an optimal cost error of $0.024\%$ (100 times more accurate than MS) and input control trajectories RMSE of $4.2\text{e-}3$, as seen in Table~\ref{tab:toy_example}. Moreover for all methods, the ODE is solved with an RMSE of $5.1\text{e-}4\%$, verifying that the continuous time OCP is solved to optimality. Computationally, the method is quicker than PS and MS, for a simple system. Notice that the PS method is able to find a smaller optimal cost, by violating the constraints between collocation points as in Figure~\ref{fig:dummy_ex_states_control}. Moreover, for higher degree splines, without a safety envelope, the solution experiences an oscillatory behavior due to over-fitting as the coefficients are not constrained intuitively. As the safety envelope is parameterized by the coefficients $\alpha$, indirect regularization is added on them.
\section{Real-time NMPC for Autonomous Valet Parking}
\label{section:5}
The presented SOCSE method was developed to speed-up and allow a real-time deployment of the NMPC for complex autonomous systems. For this reason, we will employ the method in an AVP example in which the proposed method solves considerably quicker than standard shooting and PS approaches. The AVP task is validated first in simulation then transferred to the real-world.
\subsection{Path following and autonomous valet parking control}
We represent the car dynamics in the NMPC by a 3 DoF dynamic single-track model, in the curvilinear error frame:%
{\small\begin{align}
		\label{eq:BicycleModelStates}
		\Dot{v}_x &= (F_{xf} \cos\delta + F_{xr} -F_{yf} \sin\delta - F_{res} + M\Dot{\psi} v_y)/M, \nonumber\\
		\Dot{v}_y &= (F_{xf} \sin\delta + F_{yr} + F_{yf} \cos\delta - M\Dot{\psi} v_x)/M, \\\nonumber
		\Dot{r} &= (L_f (F_{yf} \cos\delta + F_{xf} \sin\delta) - L_r F_{yr})/I_z, \\\nonumber
		\Dot{s} &= (v_x \cos\theta -v_y \sin\theta)/(1-\kappa_c w), \\\nonumber
		\Dot{w} &= v_x \sin\theta + v_y \cos\theta, \\\nonumber
		\Dot{\theta} &= \Dot{\psi} - \Dot{\psi_c} = \Dot{\psi} - \kappa_c \Dot{s},\\\nonumber
		\Dot{w}_p &= v_x \sin\theta_p + v_y \cos\theta_p, \\\nonumber
		\Dot{\theta_p} &= \Dot{\psi}. \nonumber
\end{align}}%
The first three equations of~\eqref{eq:BicycleModelStates} dictate the dynamics in the car body frame with $x$ pointing forward, $y$ to the left, and the last five dictate the kinematics in the curvilinear frame. Morevoer, $M$ is the vehicle's mass, $I_z$ is the inertia about the z-axis, and $L_r,L_f$ are the distances from the center of gravity to the rear and front axles. We augment with the states $s$ to track the evolution along the center-line, $w_p$ and $\theta_p$ for parking spot position and heading tracking errors. Moreover, $w = (Y-Y_c)\cos(\psi_c) - (X-X_c)\sin(\psi_c)$ and $\theta =  \psi - \psi_c$ are respectively the distance and heading deviation from the center-line position $(X_c, Y_c)$ and heading $\psi_c$. For smooth driving, we augment the dynamics with the input rates $u =[\dot{t}_r,\Dot{\delta}]$, the derivatives of the body frame steering angle $\delta$ and the normalized acceleration $t_r$ from front/rear axles longitudinal forces $F_{xf}$ and $F_{xr}$. We assume a linear tire model for the lateral forces $F_{yf}$ and $F_{yr}$. For control purposes, we formulate the kinematics in the curvilinear reference frame, attached to the car body frame, equivalent to an error frame with respect to the desired path. This transformation renders the NMPC formulation parameterized only by the road curvature $\kappa_c$, allowing us to constrain the deviation error in a convex form: $w_l\leq w \leq w_r$ where $w_l$ and $w_r$ are the left and right track limits.
%
Hence, the single-track curvilinear dynamics between the state vector $x = [v_x, v_y, r, s, w, \theta, w_p, \theta_p, \delta, t_r]$ and input $u$ are:%
\begin{equation}%
	\label{eq:dynamics}
		\Dot{x} = f(x,u) = \lambda f_{dyn}(x,u) + (1-\lambda)f_{kin}(x,u).
\end{equation}%
In~\eqref{eq:dynamics}, we fuse the dynamic model in~\eqref{eq:BicycleModelStates} with the kinematic one as in the former, tire slip angle and lateral force equations are not valid at low speed, leading to a numerical problem. By fusing through $\lambda(v_x)$ from~\cite{9653681}, with a switching velocity range between $0$ and $1\textrm{ms}^{-1}$, we allow driving from and to rest (zero velocity) and provide the NLP with a differentiable fused model. %
Finally, the NLP follows~\eqref{eq:SOSCSE_nlp} with the stage cost:%
{\small\begin{align}\label{eq:stage_cost_nmpc}%
		l(\tau) = (x(\tau) - x^r(\tau))^\top Q (x(\tau) - x^r(\tau))
		+ u(\tau)^\top R u(\tau) \nonumber \\ +(1-\phi)(Q_{w_p}w_p^2+Q_{\theta_p}\theta_p^2), 
\end{align}}%
with $Q = [3,10,10,0,6,30,0,0,9,1]\in\mathbb{R}^{10\times 10} \succeq 0, R = [1,1] \in\mathbb{R}^{2\times 2} \succ 0$. The path following problem is a regulation about a zero reference $x^r$ for all states except the velocity.
The parking position tracking is fused in the objective function based on the remaining distance before the parking spot $s_p$ with $\phi = \exp(-A\times(s-s_p)^2)$, $A=0.05, Q_{w_p} = 6, Q_{\theta_p} = 30$.

In simulation, we solve one step of the AVP problem, with the car initialized at the limit of the OCP boundaries, at $w = 2.99m$ from the reference center-line.  The analytical solution is approximated by solving the problem over the same horizon length $t_f=2$s, in a classical MS approach on a coarse grid of $N^*=10000$ nodes, and a RK4 integrator. To verify the solution of the ODE, we roll out the optimal control action solution $u^*(\tau)$ through a dynamics integrator with a small step size (0.1ms), that is $\tilde{x}_{k+1} = f_d(\tilde{x}_k,u^*_k)$. We then compare the resulting state $\tilde{x}$ with the optimal state trajectory solved for in the OCP $x^*(\tau)$ and provide the ODE solution error entry in Table~\ref{tab:avp_example}. %
\begin{table}[htbp]%
	\caption{Comparison of different methods for solving a complex AVP OCP}
	\begin{center}
		{\renewcommand{\arraystretch}{1.1}
			\begin{tabular}{p{0.2\linewidth} p{0.2\linewidth} p{0.2\linewidth} p{0.2\linewidth}}
				\hline
				Method & Solve time (s)& Optimal cost deviation (\%)& ODE solution error\\
				\hline
				MS-50 steps& 0.876&1.978 & - \\
				\hline
				SS-50 steps& 5.668& 1.978& - \\
				\hline
				PS $\mathcal{O}$(5)& 0.242&0.37054 & 3.384e-05 \\
				\hline
				\textbf{SOCSE $\mathcal{O}$(5)}& \textbf{0.102}&\textbf{-0.009} &  \textbf{5.978e-04}\\
				\hline
		\end{tabular}}
		\label{tab:avp_example}
	\end{center}
\end{table}%
Both multiple (MS) and single shooting (SS) methods are indistinctive from each other. However, for this particular example, SS is slower to converge than MS. Similar to~\ref{subsect:analyt_ex}, only adding inequality constraints on the collocation points results in violation of the control limits between the collocation points, and particularly at the sample time of the feedback control (50ms). However, with our proposed approach of augmenting with the safety envelope (SOCSE), the full trajectory is safe, and the optimal cost relative error is only $-0.009\%$, the smallest by a margin out of all the proposed approaches. The computation time to solve such an OCP is also the smallest at 0.1s (solved with IPOPT). Our method outperforms the standard collocation method with varying degrees, in terms of accuracy (state and control deviation with respect to MS-10000) and convergence speed, while generating a continuous safe solution as in Figure~\ref{fig:cost_states_evo}. Beyond degree 5 the solution stops improving, and higher degree coefficients are reduced to zero, given the orthogonal property of the basis.
Finally, Figure~\ref{fig:mil_control_wEnvelope} presents the results of a closed-loop simulation with a high-fidelity simulator. The shaded area shows the estimated bounds of the open-loop MPC solution at every iteration: while a standard collocation approach violates constraints between collocation points, the SOCSE method adds safety without altering the closed-loop MPC policy. 

\begin{figure}
	\begin{center}
		{
%
%
\begin{tikzpicture}

\begin{axis}[%
width=6.8cm,
height=0.878cm,
at={(0cm,3.315cm)},
scale only axis,
xmin=3,
xmax=10,
xlabel style={font=\color{white!15!black}},
ymin=-1.5,
ymax=2,
ylabel style={font=\color{white!15!black}},
ylabel style={font=\footnotesize, align=center},
ylabel={$\log_{10}$\\(Computation time)},
axis background/.style={fill=white},
axis x line*=bottom,
axis y line*=left,
xmajorgrids,
ymajorgrids
]

\addplot [color=black, line width=1.0pt, mark=asterisk, mark options={solid, black}, forget plot]
table[row sep=crcr]{%
	5	-0.5335\\
	6	0.7164\\
	7	0.1174\\
	8	0.0902\\
	9	0.3872\\
	10	1.6371\\
};
\addplot [color=blue, line width=1.0pt, mark=asterisk, mark options={solid, blue}, forget plot]
table[row sep=crcr]{%
	3	 -1.2094\\
	4	-1.1142\\
	5	-0.9711\\
	6	 -0.9879\\
	7	-0.8424\\
	8	-0.7425\\
	9	-0.5955\\
	10	-0.5615\\
};
\end{axis}

\begin{axis}[%
width=6.8cm,
height=0.878cm,
at={(0cm,1.657cm)},
scale only axis,
xmin=3,
xmax=10,
xlabel style={font=\color{white!15!black}},
ymin=-0.64275947357491,
ymax=0.343769249376924,
ylabel style={font=\footnotesize, align=center},
ylabel={Cost deviation\\(\%)},
ytick = {-0.6, 0, 0.2},
axis background/.style={fill=white},
axis x line*=bottom,
axis y line*=left,
xmajorgrids,
ymajorgrids
]
\addplot [color=black, line width=1.0pt, mark=asterisk, mark options={solid, black}, forget plot]
  table[row sep=crcr]{%
5	0.320400543939715\\
6	0.343769249376924\\
7	0.0900966644717592\\
8	0.0108598678620754\\
9	-0.0106704122942279\\
10	-0.0900739466484621\\
};
\addplot [color=blue, line width=1.0pt, mark=asterisk, mark options={solid, blue}, forget plot]
  table[row sep=crcr]{%
3	-0.64275947357491\\
4	-0.046758371424516\\
5	-0.0590228797235654\\
6	-0.0671372601763696\\
7	-0.0745878817910839\\
8	-0.0758267566855439\\
9	-0.0794236117209254\\
10	-0.0834831465768559\\
};
\end{axis}


\begin{axis}[%
width=6.8cm,
height=0.878cm,
at={(0cm,0cm)},
scale only axis,
xmin=3,
xmax=10,
xlabel style={font=\color{white!15!black}},
xlabel={Polynomial and spline degree $\mathcal{O}$},
ymin=0.0893035099547849,
ymax=0.398731301293927,
ylabel style={font=\color{white!15!black}},
ylabel style={font=\footnotesize, align=center},
ylabel={Max control\\deviation},
axis background/.style={fill=white},
axis x line*=bottom,
axis y line*=left,
xmajorgrids,
ymajorgrids
]
\addplot [color=black, line width=1.0pt, mark=asterisk, mark options={solid, black}, forget plot]
  table[row sep=crcr]{%
5	0.259540024744117\\
6	0.353115381307972\\
7	0.179759865685501\\
8	0.202424918628642\\
9	0.157959216843475\\
10	0.0893035099547849\\
};
\addplot [color=blue, line width=1.0pt, mark=asterisk, mark options={solid, blue}, forget plot]
  table[row sep=crcr]{%
3	0.398731301293927\\
4	0.185178533454216\\
5	0.174499956200978\\
6	0.174499638665714\\
7	0.174499919641941\\
8	0.17449994927339\\
9	0.174494074501583\\
10	0.171202323875344\\
};
\end{axis}


\begin{axis}[%
width=8.4cm,
height=3.3cm,
at={(-0,0)},
scale only axis,
xmin=0,
xmax=1,
ymin=0,
ymax=1,
axis line style={draw=none},
ticks=none,
axis x line*=bottom,
axis y line*=left
]
\end{axis}
\end{tikzpicture}%
}    
		\caption{AVP performance with safety envelope. The dynamics ODE and OCP are accurately solved using the SOCSE approach and scale better with higher degrees} 
		\label{fig:cost_states_evo}
	\end{center}
\end{figure}
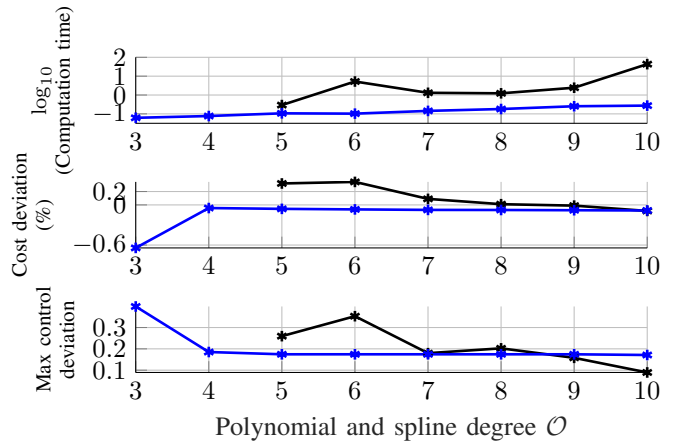
\begin{figure}
	\begin{center}
		{\input{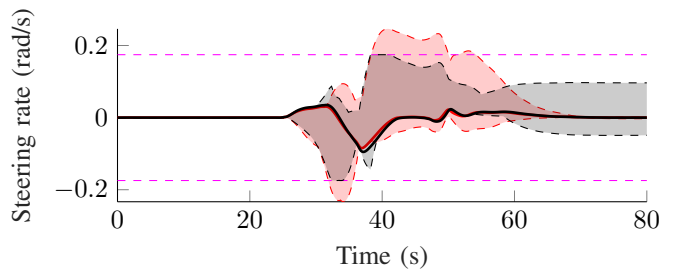}}    
		\caption{MiL AVP: Closed loop control with (black) and without (red) active safety envelope (shaded area: estimated extrema from MPC open-loop solution, pink: control limits). With a safety envelope SOCSE $\mathcal{O}$(5) (black), the full trajectory solution is contained within the control bounds} 
		\label{fig:mil_control_wEnvelope}
	\end{center}
\end{figure}


\subsection{Vehicle integration for real-time autonomous driving}
\begin{figure}
	\begin{center}
		{\input{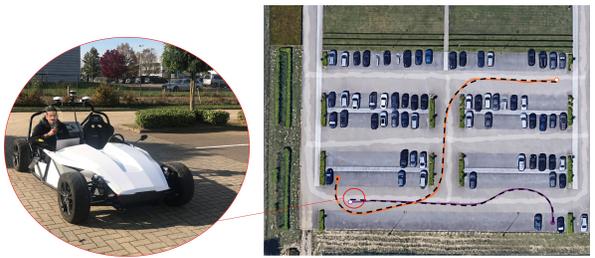}}    
		\caption{ViL AVP: tracking performance in the parking lot area for two demo runs with different paths and parking goals}
		\label{fig:vil_parking_error} 
	\end{center}
\end{figure}

The verification and validation cycle follows the methodology presented in~\cite{9867514}. We deploy our algorithm on a dSPACE MicroAutobox III which has an ARM Cortex A-15 processor, operates on a 2GB DDR4 RAM with 64MB flash memory, and runs a RT operating system. The platform in use is a drive-by-wire SimRod in Figure~\ref{fig:vil_parking_error}, equipped with a task for localization and object detection.
It is worth noting, that shooting methods were not viable for RT deployment as model-in-the-loop testing showed that they are on average 20 times slower than our approach, and did not meet the criteria to transfer to the real world. For embedded applications, the NLP is solved with an SQP method using an active-set QP solver.
In Figure~\ref{fig:vil_parking_error}, we present the results of an AVP testing campaign in a secured environment. Two trajectories are shown in the plots alongside their superposition on the real parking map. Albeit modeling errors and road grade, the NMPC was effective in both path tracking and accurate parking positioning of less than 15cm in both scenarios, with splines of degree $M=5$, using $N=6$ collocation points.
The SOCSE did not compromise performance for safety, as given the tightness of the safety envelope, the problem formulation is not conservative. 
With the SOCSE approach, we were able to 1) deploy a RT NMPC on the vehicle with a prediction horizon of 2 seconds and a sampling frequency of 20Hz 2) obtain a satisfying performance 3) autonomously drive smoothly and safely as states and inputs are smooth by construction 4) formulate the NLP independently of any sampling time.   

\section{Conclusion}
This paper presents a novel approach to formulate an NLP for optimal control problems based on collocation schemes. Without adding complexity and with minimal conservatism, the problem approximates and solves the original continuous time OCP with high accuracy, fast rate of convergence and a small number of optimization variables, all while resulting in a safe solution over the full optimization horizon. We make use of Bernstein polynomial bound approximation to efficiently transcribe the inequalities using a convex safety envelope in a spectral orthogonal collocation framework. Importantly, no discretization is needed, and the solution is scalable to any sampling time allowing a more favorable and orders faster real-time NMPC deployment than standard approaches and contributing as a first safety layer to autonomous driving.  

\bibliographystyle{IEEEtran}
\bibliography{ecc23bib}

\end{document}